\documentclass[12pt]{article}
\usepackage{amsthm}
\usepackage{index}
\usepackage[cp1251]{inputenc}
\usepackage[english]{babel}
\usepackage{amsmath, amssymb, eucal}
\usepackage{cite}
\begin{document}
	
\begin{center}
	{\bf EXISTENCE AND EXPLICIT FORM OF NONLINEAR HERMITE--CHEBYSHEV APPROXIMATIONS }
\end{center}

\begin{center}
{\bf  A.P. Starovoitov, I.V. Kruglikov}

{\it F. Skorina Gomel State University, Gomel }
\end{center}

In this paper, sufficient conditions for the existence of trigonometric Hermite--Jacobi appro\-ximations of a system of functions that are sums of convergent Fourier series are found. Based on these results, sufficient conditions are established under which nonlinear Hermite--Chebyshev approximations of systems of functions representable by Fourier series in Chebyshev polynomials of the first and second kind exist. When the found conditions are met, explicit formulas are obtained for the numerators and denominators of trigonometric Hermite--Jacobi approximations and nonlinear Hermite--Chebyshev approximations of the first and second kind of the specified systems of functions.\\


{\bf\itshape Introduction}

Let the set ${\bf f^{ch1}}=(f^{ch1}_1,\ldots,f^{ch1}_k)$ consist of functions represented by Fourier series with respect to Chebyshev polynomials $T_n(x)=\cos(n\arccos x)$ of the first kind
\begin{equation}\label{eq1}
f^{ch1}_j(x)=\frac{a_0^j}{2}+
\sum_{l=1}^{\infty} a^j_l T_l(x),\,\, j=1,\ldots,k\,,
\end{equation}
with real coefficients, that converge for all $x\in  [-1,1]$. The set of $k$--dimensional
multi-indices, which are an ordered set of $k$ non-negative integers, is denoted by $\mathbb{Z}^k_+$. The order of the multi-index $\overrightarrow{m}=(m_1,\ldots,m_k)\in \mathbb{Z}^k_+$ is the sum $m=m_1+\ldots+m_k$.

Let us fix an index $n\in
\mathbb{Z}^1_+$ and a multi-index ${\overrightarrow{m}=(m_1,
	\ldots, m_k)\in\mathbb{Z}^k_+ }$ and consider the following analogue of the Hermite--Pad\'e problem for  ${\bf f^{ch1}}$ (see \cite[chapter~4, \S 1, problem A]{NikSor}):

\vspace{0.1 cm}
{\bf Problem ${\bf A^{ch1}}$}. {\it For a system of functions ${\bf f^{ch1}}$ find a polynomial
	$
	Q^{ch1}_m(x;{\bf f^{ch1}})=Q^{ch1}_{n,\overrightarrow{m}}(x;{\bf f^{ch1}})=\sum_{p=0}^{m}u_pT_p(x)
	$ 
	that is not equal to zero identically and polynomials
	$
	P^{ch1}_j(x;{\bf f^{ch1}})=P^{ch1}_{n_j,n,\overrightarrow m}(x;{\bf f^{ch1}})=\sum_{p=0}^{n_j}v^j_p\,T_p(x)
	$, $n_j=n+m-m_j$, 
	that for $j=1,\ldots,k$
$$
Q^{ch1}_m(x;{\bf f^{ch1}})f^{ch1}_j(x)-P^{ch1}_j(x;{\bf f^{ch1}})
=\sum_{l=n+m+1}^{\infty}\tilde{a}^j_l\,T_l(x)\,.
$$}

\vspace{0,1cm} 
{\bf Definition 1.} {\it If the pair $(Q^{ch1}_m,P^{ch1})$, where $P^{ch1}=(P^{ch1}_1,\ldots,P^{ch1}_k)$, is a solution to problem ${\bf A^{ch1}}$, then rational fractions
	$$
\pi^{ch1}_j(x;{\bf f^{ch1}})=\pi^{ch1}_{n_j,n,\overrightarrow{m}}(x;{\bf f^{ch1}})=\frac{P^{ch1}_j(x;{\bf f^{ch1}})}{Q^{ch1}_m(x;{\bf f^{ch1}})},\,\,\,
	j=1,\ldots,k\,,
	$$
	will be called {\it linear Hermite--Chebyshev approximations of the first kind} for the multi-index $(n,\overrightarrow{m})$ and the system ${\bf f^{ch1}}$.}
 
\vspace{0,2cm} 
{\bf Definition 2.} {\it Nonlinear Hermite--Chebyshev approximations of the first kind for the multi-index $(n,\overrightarrow{m})$ and the system ${\bf f^{ch1}}$ will be called rational fractions
$$
\widehat{\pi}^{ch1}_j(x;{\bf f^{ch1}})=	\widehat{\pi}^{ch1}_{n_j,n,\overrightarrow{m}}(x;{\bf f^{ch1}})=
\frac{\widehat{P}^{ch1}_{j}(x;{\bf f^{ch1}})}{\widehat{Q}^{ch1}_m(x;{\bf f^{ch1}})},
$$
where the polynomials
$\widehat{Q}^{ch1}_{m}(x;{\bf f^{ch1}})$,\,\,
$\widehat{P}^{ch1}_{j}(x;{\bf f^{ch1}})$\,\, $(n_j=n+m-m_j)$, the degrees of which do not exceed $m$ and $n_j$ respectively, are chosen so that
$$
f^{ch1}_j(x)-\frac{\widehat{P}^{ch1}_{j}(x;{\bf f^{ch1}})}{\widehat{Q}^{ch1}_m(x;{\bf f^{ch1}})}=
\sum_{l=n+m+1}^{\infty} \widehat{a}^j_l\,T_l(x),\,\,\,   j=1,\ldots,k\,.
$$}

In the case when $k=1$ the basic properties of linear and nonlinear Hermite--Chebyshev approximations (in this case they are called linear and nonlinear Pad\'e--Chebyshev approximations; for more details on terminology, see \cite{Suetin1}) are described in sufficient detail (see \cite{Suetin1, Suetin2, Bek} and the literature cited there, as well as \cite{Star-Lab}--\cite{Mason2}). For example, it is well known that the linear Pad\'e--Chebyshev approximation always exists, but unlike the classical Pad\'e approximation of a power series, in general, it is not unique. The nonlinear Pad\'e--Chebyshev approximation does not always exist, but if it exists, it is always unique. Similar properties are also valid for Hermite-Chebyshev approximations~\cite{StarKechOsn,StarKrugOsn}. For $k\geqslant 1$ there are examples of systems of functions ${\bf f^{ch1}}$, for which nonlinear Hermite--Chebyshev approximations exist, but are not linear Hermite--Chebyshev approximations (see \cite{Suetin1, StarKechOsn, StarKrugOsn, Suetin3, Suetin4}). 

 \vspace{0,2cm}
  Let us now consider another type of Hermite--Chebyshev approximations. Let the set ${\bf f^{ch2}}=(f^{ch2}_1,\ldots,f^{ch2}_k)$ consist of functions represented by Fourier series with respect to Chebyshev polynomials $U_n(x)=\dfrac{1}{\sqrt{1-x^2}}\sin(n\arccos x)$ of the second kind
\begin{equation}\label{eq2}
f^{ch2}_j(x)=
\sum_{l=1}^{\infty} b^j_l U_l(x),\,\, j=1,\ldots,k,
\end{equation}
with real coefficients, that converge for all $x\in  [-1,1]$. If instead of series \eqref{eq1} we take series \eqref{eq2}, then constructions similar to the previous ones lead to linear and nonlinear Hermite--Chebyshev approximations of the second kind. An analogue of the Hermite--Pad\'e problem for series \eqref{eq2} has the form:

\vspace{0.1 cm}
{\bf Problem ${\bf A^{ch2}}$}. {\it Find a polynomial 
$Q^{ch2}_m(x;{\bf f^{ch2}})=Q^{ch2}_{n,\overrightarrow{m}}(x;{\bf f^{ch2}})$, $\deg Q^{ch2}_m\leqslant m$ that is not equal to zero identically and polynomials 
$P^{ch2}_j(x;{\bf f^{ch2}})=P^{ch2}_{n_j,n,\overrightarrow m}(x;{\bf f^{ch2}})$, $\deg P^{ch2}_j \leqslant n_j$, $n_j=n+m-m_j$, 
	that for $j=1,\ldots,k$
$$
Q^{ch2}_m(x;{\bf f^{ch2}})f^{ch2}_j(x)-P^{ch2}_j(x;{\bf f^{ch2}})	=\sum_{l=n+m+1}^{\infty}\tilde{b}^j_l\,U_l(x)\,.
$$}

\vspace{0,2cm} 
{\bf Definition 3.} {\it If the pair $(Q^{ch2}_m,P^{ch2})$, where $P^{ch2}=(P^{ch2}_1,\ldots,P^{ch2}_k)$, is a solution to the problem ${\bf A^{ch2}}$, then rational fractions
	$$
	\pi^{ch2}_j(x;{\bf f^{ch2}})=\pi^{ch2}_{n_j,n,\overrightarrow{m}}(x;{\bf f^{ch2}})=\frac{P^{ch2}_j(x;{\bf f^{ch2}})}{Q^{ch2}_m(x;{\bf f^{ch2}})},\,\,\,
	j=1,\ldots,k\,,
	$$
	will be called {\it linear Hermite--Chebyshev approximations of the second kind} for the multi-index $(n,\overrightarrow{m})$ and the system ${\bf f^{ch2}}$.}

\vspace{0,2cm} 
{\bf Definition 4.} {\it Nonlinear Hermite--Chebyshev approximations of the second kind for the multi-index $(n,\overrightarrow{m})$ and the system ${\bf f^{ch2}}$ will be called algebraic rational fractions
	$$
	\widehat{\pi}^{ch2}_j(x;{\bf f^{ch2}})=	\widehat{\pi}^{ch2}_{n_j,n,\overrightarrow{m}}(x;{\bf f^{ch2}})=
	\frac{\widehat{P}^{ch2}_{j}(x;{\bf f^{ch2}})}{\widehat{Q}^{ch2}_m(x;{\bf f^{ch2}})},
	$$
	where the polynomials
	$\widehat{Q}^{ch2}_{m}(x;{\bf f^{ch2}})$,\,\,
	$\widehat{P}^{ch2}_{j}(x;{\bf f^{ch2}})$\,\, $(n_j=n+m-m_j)$, the degrees of which do not exceed respectively $m$ and $n_j$, are chosen so that
	$$
	f^{ch2}_j(x)-\frac{\widehat{P}^{ch2}_{j}(x;{\bf f^{ch2}})}{\widehat{Q}^{ch2}_m(x;{\bf f^{ch2}})}=
	\sum_{l=n+m+1}^{\infty} \widehat{b}^j_l\,U_l(x),\,\,\,  j=1,\ldots,k\,.
	$$}

Now consider two systems ${\bf f^{1}}=(f^{1}_1,\ldots,f^{1}_k)$, ${\bf f^{t1}}=(f^{t1}_1,\ldots,f^{t1}_k)$ of power and trigonometric series
$$
f^{1}_j(x)=\frac{a_0^j}{2}+
\sum_{l=1}^{\infty} a^j_l z^l,\,\,\,
f^{t1}_j(x)=\frac{a_0^j}{2}+
\sum_{l=1}^{\infty} a^j_l \cos lx, 
$$
{\it associated with the system ${\bf f^{ch1}}$},
and two systems ${\bf f^{2}}=(f^{2}_1,\ldots,f^{2}_k)$, ${\bf f^{t2}}=(f^{t2}_1,\ldots,f^{t2}_k)$ of power and trigonometric series
$$
f^{2}_j(x)=\sum_{l=1}^{\infty} b^j_l z^l,\,\,\,
f^{t2}_j(x)=\sum_{l=1}^{\infty} b^j_l \sin lx,
$$
{\it associated with the system ${\bf f^{ch2}}$}.

\vspace{0,2cm}
 For $k=1$ in~\cite{Suetin1}\,, relying on the properties of Faber operators, it is established that the question of the existence of nonlinear Pad\'e--Chebyshev approximations of the first kind is resolved using known results on classical Pad\'e approximations of the corresponding associated power series. The results obtained in this paper allow us to conclude that the question of the existence of nonlinear Hermite--Chebyshev approximations of the first kind can be similarly resolved using known results on Hermite--Pad\'e approximations of the corresponding associated system of power series. The main theorems of the paper were obtained by us without using the theory of Faber operators, and the proposed method for solving the problem turned out to be applicable to Hermite--Chebyshev approximations of the second kind.
 If the proof of the main result in~\cite{Suetin1} is illustrated by the scheme: 
 {\it Pad\'e approximations of a power series $\Longleftrightarrow$ Faber transformation $\Longleftrightarrow$ Pad\'e--Chebyshev approximations}, then the scheme of our reasoning looks like this:
{\it Hermite--Pad\'e approximations of a system of power series $\Longleftrightarrow$ trigonometric Hermite--Jacobi approximations $\Longleftrightarrow$ Hermite--Chebyshev approximations}. The role of Faber transformations in our scheme of reasoning is occupied by trigonometric Hermite--Jacobi approximations, to which a separate paragraph is devoted.

\vspace{0,1cm} 
Further, we will consider only {\it nonlinear} Hermite--Chebyshev approximations, and the main objective of the research in this paper is to find conditions on the coefficients of the series \eqref{eq1} and \eqref{eq2} under which nonlinear Hermite--Chebyshev approximations of the first and second kind exist. If they exist, we will look for an explicit form of such approximations. The proof of the main theorems of the paper relies heavily on the connection established in \cite{StarKechOsn, StarKechOsn23, StarKrugOsn} between nonlinear Hermite--Chebyshev approximations and the corresponding trigonometric Hermite--Jacobi approximations.

 Another of our papers \cite{StarKrug25} is devoted to the description of the conditions for the existence and uniqueness of linear Hermite--Chebyshev approximations of the first and second kind.

\vspace{0,2cm}
{\bf\itshape 1. Hermite--Pad\'e and Hermite--Jacobi approximations }

Let us present some well-known facts from the theory of Hermite--Pad\'e and Hermite--Jacobi approximations, which will be needed later.

Let ${\bf f}=(f_1,\ldots,f_k)$ be a set of $k$ power series
\begin{equation}\label{eq3}
f_j(z)=
\sum_{l=0}^{\infty}{f^j_l}{z^{l}},\,\,\,
                              j=1,\ldots,k\,,
\end{equation}
with complex coefficients. 
Let us fix $n\in
\mathbb{Z}^1_+$ and ${\overrightarrow{m}=(m_1,
	\ldots, m_k)\in\mathbb{Z}^k_+ }$ and consider the well-known Hermite--Pad\'e problem~\cite[chapter~4, \S 1, problem A]{NikSor}:

\vspace{0,2cm} {\bf Problem A}. {\it Find a polynomial
	$Q_m(z;{\bf f})=Q_{n,\overrightarrow{m}}(z;{\bf f})$ that is not identically equal to zero and whose degree $\deg Q_m \leqslant m$, and polynomials $P_{j}(z;{\bf f})=P_{n_j,n,\overrightarrow m}(z;{\bf f})$, 
	$\deg P_j\leqslant n_j$, $n_j=n+m-m_j$, so that for}
$j=1,\ldots,k$
\begin{equation}\label{eq4}
Q_m(z;{\bf f})f_j(z)-P_j(z;{\bf f})=O(z^{n+m+1}).
\end{equation}
From here on, by $O(z^p)$ we mean a power series of the form $c_1z^p+c_2z^{p+1}+\ldots$\,\,\,.

\vspace{0,2cm} 
{\bf Definition 5.} {\it If the polynomials	$Q_m(z;{\bf f}), P_1(z;{\bf f}),\ldots, P_k(z;{\bf f})$ are solutions to problem~{\bf A} (a solution to problem~{\bf A} always exists~\cite{NikSor}), 
then the rational fractions
$$
\pi_j(z;{\bf f})=	\pi_{n_j,n,\overrightarrow{m}}(z;{\bf f})=\frac{P_j(z;{\bf f})}{Q_m(z;{\bf f})},\,\,\,j=1,\ldots,k\,,
$$
are called {\it Hermite--Pad\'e approximations} for the multi-index $(n,\overrightarrow{m})$ and the system ${\bf f}$.}

 For the system of exponential functions, the fractions $\{\pi_j(z;{\bf f})\}_{j=1}^{k}$ were first introduced into consideration by Ch.\,Hermite in his paper~\cite{Hermite}, devoted to proving the transcendence of the number $e$. Conditions \eqref{eq4} do not uniquely determine them \cite{NikSor, StarRjab2}. 
 Of particular interest are systems of functions $\bf f$, for which
 $\{\pi_j(z;{\bf f})\}_{j=1}^k$ from \eqref{eq4} are determined uniquely up to a normalizing factor for any multi-index $(n,\overrightarrow{m})$. Such systems are, for example, perfect systems.~\cite[chapter~4, \S 1]{NikSor}. There are systems ${\bf f}$, other than perfect ones, for which $\{\pi_j(z;{\bf f})\}_{j=1}^k$ are also uniquely determined~\cite{StarRjab1, StarRjab2}.

\vspace{0,2cm}
Let us introduce into consideration multiple analogues of the fractions of C.\,Jacobi~\cite{Jacobi} (for more details see~\cite{ApSuet}).

\vspace{0,1cm} 
{\bf Definition 6.} {\it
	Rational functions of the form
$$
\widehat{\pi}_j(z;{\bf f})=	\widehat{\pi}_{n_j,n,\overrightarrow{m}}(z;{\bf f})=\frac{\widehat{P}_j(z;{\bf f})}{\widehat{Q}_m(z;{\bf f})},\,\,\,
j=1,\ldots,k\,,
$$
	where algebraic polynomials $\widehat{Q}_m(z;{\bf f})=\widehat{Q}_{n,\overrightarrow{m}}(z;{\bf f})$,
	$\widehat{P}_j(z;{\bf f})=\widehat{P}_{n_j,n,\overrightarrow{m}}^j(z;{\bf f})$ have degrees not higher than $m$ and $n_j$, $n_j=n+m-m_j$ respectively, will be called {\it Hermite--Jacobi approximations} for a multi-index $(n,\overrightarrow{m})$ and a system of functions ${\bf f}$, if 
$$
	f_j(z)-\frac{\widehat{P}_j(z;{\bf f})}{\widehat{Q}_m(z;{\bf f})}=O(z^{n+m+1}).
$$}

Unlike Hermite--Pad\'e approximations, Hermite--Jacobi approximations may not exist \cite{StarKrugOsn}.

Let us introduce new notations. For index $n\in
\mathbb{Z}^1_+$ multi-index ${\overrightarrow{m}=(m_1,
	\ldots, m_k) }$, $m \neq 0$, consider a determinant

\begin{equation*}
H_{n,\overrightarrow m}({\bf f})=\det
\left(  \begin{array}{cccc}
H^1 \\
H^2 \\
\vdots  \\
H^k \\
\end{array}  \right)
=
\left|  \begin{array}{cccc}
f_{n-m_1+1}^{1} & f_{n-m_1+2}^{1} & \ldots & f_{n_1}^{1} \\
\ldots & \ldots & \ldots & \ldots \\
f_{n}^{1} & f_{n+1}^{1} & \ldots & f_{n+m-1}^{1} \\
\ldots & \ldots & \ldots & \ldots \\
f_{n-m_k+1}^k & f_{n-m_k+2}^k & \ldots & f_{n_k}^k \\
\ldots & \ldots & \ldots & \ldots \\
f_{n}^k & f_{n+1}^k & \ldots & f_{n+m-1}^k \\
\end{array}  \right|,
\end{equation*}
for  $m_j\neq0$ consisting of blocks
$$
H^j=\left( \! \begin{array}{cccc}
f_{n-m_j+1}^{j} & f_{n-m_j+2}^{j} & \ldots & f_{n_j}^{j} \\
f_{n-m_j+2}^{j} & f_{n-m_j+3}^{j} & \ldots & f_{n_j+1}^{j} \\
\ldots & \ldots & \ldots & \ldots \\
f_{n}^j & f_{n+1}^j & \ldots & f_{n+m-1}^j \\
\end{array} \! \right),
$$
located one above the other.  For $p<0$ we assume that $f^j_p=0$. By definition, we assume that for $m_j=0$, determinant $H_{n,\overrightarrow m}({\bf f})$ does not contain block $H^j$. For $k=1$ the determinant $H_{n,\overrightarrow m}({\bf f})$ is the well-known Hadamard determinant (see \cite{Bek}).

In \cite{OsnRyabStar} (see also \cite{StarKrugOsn}) a multiple analogue of Jacobi's theorem \cite{Jacobi, Bek} was proved: if for a multi-index $(n,\overrightarrow{m})$ and a system of functions ${\bf f}$ the determinant $H_{n,\overrightarrow{m}}({\bf f}) \neq 0$, then the Hermite--Jacobi approximations $\{\widehat{\pi}_j(z;{\bf f})\}_{j=1}^k$ exist, are uniquely determined, and each of them identically coincides with the corresponding Hermite--Pad\'e approximation, i.e.
$$
\widehat{\pi}_j(z;{\bf f})=\pi_j(z;{\bf f}),\,\,\, j=1,\ldots,k.
$$
 In this case, based on the results of paper~\cite{StarRjab2}, we can describe the explicit form of the numerators and denominators of the fractions $\{\widehat{\pi}_j(z;{\bf f})\}_{j=1}^k$. For this, we introduce the necessary notation.
 
  For  $m_j\neq 0$, we add the column
 $$
 \left( \! \begin{array}{cccccc}
 f_{n_j+1}^{j} & f_{n_j+2}^{j} & \ldots  & f_{n+m}^{j}\\
 \end{array} \! \right)^T
 $$
 as the last column of the matrix $H^j$.
As a result, we obtain a matrix $F^j$ of order $m_j\times(m+1)$. By placing $F^j$ one above the other according to their numbers, we construct a new matrix of order $m\times(m+1)$
$$
F_{n,\,\overrightarrow{m}}({\bf f})=\left[ \! \begin{array}{cccc}
F^1 & F^2 & \ldots & F^k\\
\end{array} \!\, \right]^{T}\,:=\left( \! \begin{array}{cccc}
f_{n-m_1+1}^{1} & f_{n-m_1+2}^{1} & \ldots & f_{n_1+1}^{1} \\
f_{n-m_1+2}^{1} & f_{n-m_1+3}^{1} & \ldots & f_{n_1+2}^{1} \\
\ldots & \ldots & \ldots & \ldots \\
f_{n}^{1} & f_{n+1}^{1} & \ldots & f_{n+m}^{1} \\
\ldots & \ldots & \ldots & \ldots \\
f_{n-m_k+1}^k & f_{n-m_k+2}^k & \ldots & f_{n_k+1}^k \\
f_{n-m_k+2}^k & f_{n-m_k+3}^k & \ldots & f_{n_k+2}^k \\
\ldots & \ldots & \ldots & \ldots \\
f_{n}^k & f_{n+1}^k & \ldots & f_{n+m}^k \\
\end{array} \! \right).                       
$$
We also define functional row matrices of order
$1\times(m+1)$
\begin{equation*}
E(z)=\left(  \begin{array}{cccccc}
z^m & z^{m-1} & \ldots & z & 1\\
\end{array}  \right)\,,
\end{equation*}
\begin{equation*}
E_{m_j}(z)=\left(  \begin{array}{cccccc}
\sum\limits_{l=0}^{n-m_j}f_l^{j}z^{m+l} &
\sum\limits_{l=0}^{n-m_j+1}f_l^{j}z^{m+l-1} & \ldots & \sum\limits_{l=0}^{n_j}f_l^{j}z^{l}\\
\end{array}  \right)\,.
\end{equation*}
If we add the rows $E(z)$, $E_{m_j}(z)$ and
$
\left( \! \begin{array}{cccccc}
	f_{n+l}^{j} & f_{n+l+1}^{j} & \ldots  & f_{n+m+l}^{j}\\
\end{array} \! \right),
$
to the matrix $F_{n,\,\overrightarrow{m}}({\bf f})$ as the last row, respectively, we obtain square matrices whose determinants we denote by $D(n,\overrightarrow{m};z)$, $G_j(n,\overrightarrow{m};z)$ and $d_{n,
	\overrightarrow{m},l}^j$\,, respectively.

If $H_{n,\overrightarrow{m}}({\bf f}) \neq 0$, then $F_{n,\,\overrightarrow{m}}({\bf f})$ is a full rank matrix and then, with a suitable choice of the normalizing factor, the denominator and numerator of the fractions $\pi_j(z;{\bf f})$ and $\widehat{\pi}_j(z;{\bf f})$ for $j=1,\ldots,k$ can be represented in the form~\cite{StarRjab2}:
\begin{equation}\label{eq5}
{Q}_m(z;{\bf f})=\widehat{Q}_m(z;{\bf f})=D(n,\overrightarrow{m};z)=\det\left[ \!\!\! \begin{array}{ccccc}\,\, F^{1} & F^{2} &
\ldots & F^k &E(z)\ \
\end{array} \!\!\!\, \right]^{T}\,\,,
\end{equation}
\begin{equation}\label{eq6}
{P}_j(z;{\bf f})=\widehat{P}_j(z;{\bf f})=G_j(n,\overrightarrow{m};z)=\det\left[ \!\!\! \begin{array}{ccccc} \,\,F^{1} &
F^{2} & \ldots & F^k &E_{m_j}(z)\ \
\end{array} \!\!\!\, \right]^{T}\,,
\end{equation}
besides,
\begin{equation}\label{eq7}
Q_m(z;{\bf f})f_j(z)-P_j(z;{\bf f})=\sum_{l=1}^{\infty}d_{n,
	\overrightarrow{m},l}^jz^{n+m+l}\,\,.
\end{equation}
If the coefficients of the series \eqref{eq3} are real numbers, then the determinants $d_{n,\overrightarrow{m},l}^j$ take real values, and $Q_m(z;{\bf f})$, ${P}_j(z;{\bf f})$ are algebraic polynomials with real coefficients~\cite{StarRjab2}. 
From \eqref{eq7} it follows that for $j=1,\ldots,k$ 
$$
f_j(z)-\widehat{\pi}_j(z;{\bf f})=\sum_{l=1}^{\infty}\widehat{d}_{n,
	\overrightarrow{m},l}^jz^{n+m+l}\,,
$$
in which the coefficients $\widehat{d}_{n,\overrightarrow{m},l}^j$ are real numbers.

\vspace{0,2cm}
{\bf\itshape 2. Trigonometric Hermite--Pad\'e and Hermite--Jacobi approximations}

\vspace{0.2 cm}
Let ${\bf f^t}=(f^t_1,...,f^t_k)$ be a set of trigonometric series
\begin{equation}\label{eq8}
f^t_j(x)=\frac{a^j_0}{2}+\sum_{l=1}^{\infty}
\left (a^j_l\cos lx+b^j_l\sin lx\right ),\,\,\,
j=1,\ldots,k,
\end{equation}
with real coefficients. We assume that the series~\eqref{eq8} converge for all $x\in \mathbb{R}$ and each series defines a function defined on the entire real line.

For the system ${\bf f^t}$ there exists a \cite{StarKechOsn} identically non-zero trigonometric polynomial
$Q^t_m(x;{\bf f^t})=Q^t_{n,\overrightarrow{m}}(x;{\bf f^t})$, $\deg Q^t_m \leqslant m$ and
trigonometric polynomials $P^t_{j}(x;{\bf f^t})=P^t_{n_j,n,\overrightarrow m}(x;{\bf f^t})$,
$\deg P^t_{j}\leqslant n_j$, $n_j=n+m-m_j$, such that for
$j=1,...,k$
\begin{equation}\label{eq9}
Q^t_m(x;{\bf f^t})f^t_j(x)-P^t_{j}(x;{\bf f^t})
=\sum_{l=n+m+1}^{\infty}
(\tilde{a}^j_l\cos lx+\tilde{b}^j_l\sin lx ),\,\,\,
j=1,\ldots,k.
\end{equation}

\vspace{0,2cm} 
{\bf Definition 7.} {\it If the polynomials $Q^t_m(x;{\bf f^t})$, $P^t_{j}(x;{\bf f^t})$, $j=1,...,k$, satisfy conditions \eqref{eq9}, then the trigonometric rational fractions
	$$
	\pi^t_j(x;{\bf f^t})=\pi^t_{n_j,n,\overrightarrow{m}}(x;{\bf f^t})=\frac{P^t_j(x;{\bf f^t})}{Q^t_m(x;{\bf f^t})},\,\,\,
	j=1,\ldots,k\,,
	$$
	will be called {\it trigonometric Hermite--Pad\'e approximations} ({\it joint Hermite--Fourier approximations}) for the multi-index $(n,\overrightarrow{m})$ and the system ${\bf f^t}$.}

\vspace{0,2cm}
 Note that, generally speaking, the trigonometric Hermite--Pad\'e approximations are not uniquely determined by conditions \eqref{eq9}. The problem of finding conditions under which they are uniquely determined is researched in detail in the papers \cite{Star-Lab, StarKechOsn23, StarKechOsn}.

\vspace{0.1 cm}
Let us now define the trigonometric analogues of the Hermite--Jacobi approximations.

\vspace{0,2cm} 
{\bf Definition 8.} {\it Trigonometric rational functions
$$
\widehat{\pi}^t_j(x;{\bf f^t})=\widehat{\pi}^t_{n_j,n,\overrightarrow{m}}(x;{\bf f^t})=\frac{\widehat{P}^t_{j}(x;{\bf f^t})}{\widehat{Q}^t_m(x;{\bf f^t})},\,\,\,
	j=1,\ldots,k\,,
$$
where $\widehat{Q}^t_m(x;{\bf f^t})=\widehat{Q}^t_{n,\overrightarrow{m}}(x;{\bf f^t})$,  
 $\widehat{P}^t_{j}(x;{\bf f^t})=\widehat{P}^t_{n_j,n,\overrightarrow m}(x;{\bf f^t})$ are trigonometric polynomials whose degrees, respectively, are not higher than $m$ and $n_j$, $n_j=n+m-m_j$,
will be called {\it trigonometric Hermite--Jacobi approximations} for the multi-index $(n,\overrightarrow{m})$ and the system ${\bf f^t}$, if 
$$
f^t_j(x)-\frac{\widehat{P}^t_{j}(x;{\bf f^t})}{\widehat{Q}^t_m(x;{\bf f^t})}=\sum_{l=n+m+1}^{\infty}
(\widehat{a}^j_l\cos lx+\widehat{b}^j_l\sin lx ),\,\,\,
j=1,\ldots,k.
$$}

\vspace{0,2cm}
Unlike trigonometric Hermite--Pad\'e approximations, trigonometric Hermite--Jacobi approximations may not exist \cite{Star-Lab, Suetin1, Suetin3}. There are examples of systems of functions ${\bf f^t}$, for which trigonometric Hermite--Jacobi approximations exist, but are not trigonometric Hermite--Pad\'e approximations \cite{StarKechOsn, StarKrugOsn, Nemeth}.

\vspace{0,2cm}
{\bf Theorem 1.} {\it Suppose that the multi-index $(n,\overrightarrow{m})$,  $m\neq 0$, satisfies the condition $n \geqslant \max \{ m_j\!: 1\leqslant j \leqslant k \}$, and
 ${\bf f^{1}}=(f^{1}_1,\ldots,f^{1}_k)$, ${\bf f^{t1}}=(f^{t1}_1,\ldots,f^{t1}_k)$ are two systems of power and trigonometric series associated with ${\bf f^{ch1}}$. 	
Then, for the existence of trigonometric Hermite--Jacobi approximations $\{\widehat{\pi}^t_j(x;{\bf f^{t1}})\}_{j=1}^k$ it is sufficient that the following conditions are satisfied for the system of power series $\bf f^1$:
	
	1) there exist algebraic Hermite--Jacobi approximations $\{\widehat{\pi}_j(z;{\bf f^1})\}_{j=1}^k$;
	
	2) each power series $f^1_j$ of the system $\bf f^1$ has a radius of convergence $R^1_j>1$;
	
	3) the rational fractions $\widehat{\pi}_j(z;{\bf f^1})$, $j=1,\ldots,k$, have no poles in $\overline{D}= \{ z : |z| \leqslant 1 \}$.
	
	If conditions 1)--3) are satisfied for the system $\bf f^1$, then with the appropriate normalization for all $j=1,\dots,k$ the following equalities are true:}
$$
\widehat{Q}^t_m(x;{\bf f^{t1}})=\widehat{Q}_m(e^{ix};{\bf f^1}) \cdot \overline{\widehat{Q}_m(e^{ix};{\bf f^1})}=\sum_{l=0}^m{q_l\cos lx}, 
$$
$$
\widehat{P}^t_j(x;{\bf f^{t1}})= \textrm{Re} \left \{ \widehat{P}_j(e^{ix};{\bf f^1}) \cdot \overline{\widehat{Q}_m(e^{ix};{\bf f^1})} \right \}=\sum_{l=0}^{n_j}{p^j_l\cos lx}, 
$$
$$
f^{t1}_j(x)-\widehat{\pi}^t_j(x;{\bf f^{t1}})=\sum_{l=n+m+1}^{\infty}{d^j_l\cos lx}, 
$$
{\it where $q_l$, $p_l^j$, $d_l^j$ are real numbers.}

\vspace{0,2cm}
{\bf Theorem 2.} {\it Suppose that the multi-index $(n,\overrightarrow{m})$,  $m\neq 0$, satisfies the condition $n\geqslant \max\{m_j\!:1\leqslant j \leqslant k \}$, and
	${\bf f^{2}}=(f^{2}_1,\ldots,f^{2}_k)$, ${\bf f^{t2}}=(f^{t2}_1,\ldots,f^{t2}_k)$ are two systems of power and trigonometric series associated with ${\bf f^{ch2}}$. 	
	Then, for the existence of trigonometric Hermite--Jacobi approximations $\{\widehat{\pi}^t_j(x;{\bf f^{t2}})\}_{j=1}^k$ it is sufficient that the following conditions are satisfied for the system of power series $\bf f^2$:
	
	1) there exist algebraic Hermite--Jacobi approximations $\{\widehat{\pi}_j(z;{\bf f^2})\}_{j=1}^k$;
	
	2) each power series $f^2_j$ of the system $\bf f^2$ has a radius of convergence $R^2_j>1$;
	
	3) the rational fractions $\widehat{\pi}_j(z;{\bf f^2})$, $j=1,\ldots,k$, have no poles in $\overline{D}= \{ z : |z| \leqslant 1 \}$.
	
	If conditions 1)--3) are satisfied for the system $\bf f^2$, then with the appropriate normalization for all $j=1,\dots,k$ the following equalities are true:}
	\begin{equation}\label{eq10}
	\widehat{Q}^t_m(x;{\bf f^{t2}})=\widehat{Q}_m(e^{ix};{\bf f^2}) \cdot \overline{\widehat{Q}_m(e^{ix};{\bf f^2})}=\sum_{l=0}^m{\widehat{q_l}\cos lx}, 
	\end{equation}
	\begin{equation}\label{eq11}
	\widehat{P}^t_j(x;{\bf f^{t2}})= \textrm{Im} \left \{ \widehat{P}_j(e^{ix};{\bf f^2}) \cdot \overline{\widehat{Q}_m(e^{ix};{\bf f^2})} \right \}=\sum_{l=0}^{n_j}{\widehat{p}^j_l\sin lx},
	\end{equation}
	\begin{equation}\label{eq12}
	f^{t2}_j(x)-\widehat{\pi}^t_j(x;{\bf f^{t2}})=\sum_{l=n+m+1}^{\infty}{\widehat{d}^j_l\sin lx},
	\end{equation}
	{\it where $\widehat{q}_l$, $\widehat{p}_l^j$, $\widehat{d}_l^j$ are real numbers.}

\vspace{0,2cm}
{\bf Remark 1.} Conditions 2) and 3) in theorems 1 and 2 are not essential in the following sense. If they are not satisfied, for example, for ${\bf f^{1}}$, then we can move on to another system of functions ${\bf f^{1}_r}=(f^{1}_1(rz),\ldots,f^{1}_k(rz))$, for which, for a sufficiently small $r$, $0<r<1$, they will be satisfied.

\vspace{0,2cm}
{\bf Remark 2.} In the formulations of theorems 1 and 2 it is assumed that the multi-index satisfies the condition $n\geqslant \max\{ m_j\!:1\leqslant j \leqslant k \}$. This means that we consider only the "upper"\, part of the general trigonometric Pad\'e table (see \cite[chapter~4, \S 1]{NikSor}), including its main diagonal. This condition on the multi-index allows us to significantly simplify the finding of explicit formulas for Hermite--Pad\'e approximations and is therefore found in a number of papers devoted to this topic (see, for example, \cite{Suetin1, Ibr, Bruin, Aptek81M-L, Star2018}).

\vspace{0,2cm}
Let us dwell only on the proof of theorem 2. Theorem 1 is proved similarly.
Since condition 1) is satisfied and the coefficients of the series \eqref{eq2} are real numbers, the fractions
$\{\widehat{\pi}_j(z;{\bf f^2})\}_{j=1}^k$ are defined, and their numerators and denominators are polynomials with real coefficients. Therefore, in some
neighborhood of zero
\begin{equation}\label{eq13}
f^2_j(z)-\widehat{\pi}_j(z;{\bf f^2})=
\sum_{l=n+m+1}^{\infty}
{\widehat{d}^j_lz^l},\,\,\,j=1,\dots,k.
\end{equation}
Satisfaction of conditions 2) and 3) allows us to take as such a neighborhood an open circle with a center at zero, the radius of which is greater than $1$. It is obvious that $f^{t2}_j(x)=\textrm{Im} \, f^2_j(e^{ix})$. Then, setting $z=e^{ix}$ in \eqref{eq13} and then equating the imaginary parts of the new equality, we obtain
\begin{equation}\label{eq14}
f^{t2}_j(x)-\textrm{Im} \left 
\{ \widehat{\pi}_j(e^{ix};{\bf f^2}) \right \} =\sum_{l=n+m+1}^{\infty}
{\widehat{d}^j_l\sin lx}.
\end{equation}
 Let us show that for $n \geqslant \max \{  m_j \!: 1 \leqslant j \leqslant k \}$
\begin{equation}\label{eq15}
\widehat{\pi}^t_j(x;{\bf f^{t2}})=\textrm{Im} 
\left \{ \widehat{\pi}_j(e^{ix};
{\bf f^2}) \right \}.
\end{equation}
 The denominator $\widehat{Q}_m(z;{\bf f^2})$ and the numerator $\widehat{P}_j(z;{\bf f^2})$ of the fraction $\widehat{\pi}_j(z;{\bf f^2})$ are algebraic polynomials with real coefficients. Let's assume that they are represented in the form
$$
\widehat{Q}_m(z;{\bf f^2})=\sum_{l=0}^{m}{\widehat{q}_lz^l}, \,\,\, \widehat{P}_j(z;{\bf f^2})=\sum_{l=0}^{n_j}{\widehat{p}^j_lz^l}.
$$
Then for $z=e^{ix}$
$$
\textrm{Im} \left \{ \widehat{\pi}_j(e^{ix};{\bf f^2}) \right \}=\frac{1}{2i} \left( \frac{\widehat{P}_j(e^{ix};{\bf f^2})}{\widehat{Q}_m(e^{ix};{\bf f^2})} - \frac{\overline{\widehat{P}_j(e^{ix};{\bf f^2})}}{\overline{\widehat{Q}_m(e^{ix};{\bf f^2})}} \right)=
$$
\begin{equation}\label{eq16}
=\frac{1}{2i}\frac{\sum_{l=0}^{n_j}{\widehat{p}^j_le^{ilx}} \cdot \sum_{s=0}^{m}{\widehat{q}_se^{-isx}} - \sum_{l=0}^{n_j}{\widehat{p}^j_le^{-ilx}} \cdot \sum_{s=0}^{m}{\widehat{q}_se^{isx}}}
{\sum_{s=0}^{m}{\widehat{q}_se^{isx}}
	 \cdot \sum_{l=0}^{m}
	{\widehat{q}_le^{-ilx}}}=
\end{equation}
$$
=\frac{\sum_{l=0}^{n_j}{\sum_{s=0}^{m}{\widehat{p}^j_l\widehat{q}_s\sin (l-s)x}}}{\sum_{s=0}^{m}{\sum_{l=0}^{m}{\widehat{q}_s\widehat{q}_l\cos (l-s)x}}}.
$$
If $n \geqslant \max \{  m_j : 1 \leqslant j \leqslant k \}$, then $n_j\geqslant m$. Therefore, the trigonometric polynomial in the numerator of fraction \eqref{eq16} has a degree not higher than $n_j$, and the trigonometric polynomial in the denominator has a degree not higher than $m$. From here and from \eqref{eq14} we conclude that equality \eqref{eq15} is valid. Then from \eqref{eq14}--\eqref{eq16} follows the validity of equalities \eqref{eq10}--\eqref{eq12}. Theorem 2 is proved.

\vspace{0,2cm}
{\bf\itshape 3. Existence of nonlinear Hermite--Chebyshev approximations}

Theorems 1 and 2 allow us to find sufficient conditions under which nonlinear Hermite--Chebyshev approximations exist and to describe their explicit form.

\vspace{0.2 cm}
{\bf Theorem 3.} \textit{Let  $n \geqslant \max\{m_j\!:1 \leqslant j \leqslant k \}$, $m \neq 0$, and for the system ${\bf f^{1}}$ of power series associated with the system ${\bf f^{ch1}}$, conditions 1)--3) of theorem 1 are satisfied. Then there exist nonlinear Hermite--Chebyshev approximations of the first kind $\{\widehat{\pi}^{ch1}_j(x;{\bf f^{ch1}})\}_{j=1}^k$ and for their denominator and numerators the following representations are valid: }
$$
\widehat{Q}^{ch1}_m(x;{\bf f^{ch1}})=\widehat{Q}_m(e^{i\arccos x};{\bf f^1}) \cdot \overline{\widehat{Q}_m(e^{i\arccos x};{\bf f^1})},
$$
$$
\widehat{P}^{ch1}_j(x;{\bf f^{ch1}})=\textrm{Re} \left\{ \widehat{P}_j(e^{i\arccos x};{\bf f^1}) \cdot \overline{\widehat{Q}_m(e^{i\arccos x};{\bf f^1})} \right\},
$$
{\it where the Hermite--Jacobi polynomials $\widehat{Q}_m(z;{\bf f^1})$,
$\widehat{P}_j(z;{\bf f^1})$ coincide with the Hermite--Pad\'e polynomials ${Q}_m(z;{\bf f^1})$,
${P}_j(z;{\bf f^1})$ and in the case when the determinant $H_{n,\overrightarrow{m}}({\bf f^1}) \neq 0$, are found by formulas \eqref{eq5}, \eqref{eq6}, in which ${\bf f}={\bf f^1}$.}

\vspace{0.2 cm}
{\bf Theorem 4.} \textit{Let  $n \geqslant \max\{m_j\!:1 \leqslant j \leqslant k \}$, $m \neq 0$, and for the system ${\bf f^{2}}$ of power series associated with the system ${\bf f^{ch2}}$, conditions 1)--3) of theorem 2 are satisfied. Then there exist nonlinear Hermite--Chebyshev approximations of the second kind $\{\widehat{\pi}^{ch2}_j(x;{\bf f^{ch2}})\}_{j=1}^k$ and for their denominator and numerators the following representations are valid: }
\begin{equation}\label{eq17}
\widehat{Q}^{ch2}_m(x;{\bf f^{ch2}})=\widehat{Q}_m(e^{i\arccos x};{\bf f^2}) \cdot \overline{\widehat{Q}_m(e^{i\arccos x};{\bf f^2})},
\end{equation}
\begin{equation}\label{eq18}
\widehat{P}^{ch2}_j(x;{\bf f^{ch2}})=\frac{1}{\sqrt{1-x^2}}\textrm{Im} \left\{ \widehat{P}_j(e^{i\arccos x};{\bf f^2}) \cdot \overline{\widehat{Q}_m(e^{i\arccos x};{\bf f^2})} \right\},
\end{equation}
{\it where the Hermite--Jacobi polynomials $\widehat{Q}_m(z;{\bf f^2})$,
$\widehat{P}_j(z;{\bf f^2})$ coincide with the Hermite--Pad\'e polynomials ${Q}_m(z;{\bf f^2})$,
${P}_j(z;{\bf f^2})$ and in the case when the determinant $H_{n,\overrightarrow{m}}({\bf f^2}) \neq 0$, are found by formulas \eqref{eq5}, \eqref{eq6}, in which ${\bf f}={\bf f^2}$.}

\textit{Proof} of theorem 4. Consider the system ${\bf f^{ch2}}$ associated with the system ${\bf f^{t2}}$. Since the conditions of theorem 2 are satisfied, there exist trigonometric Hermite--Jacobi approximations $\{\widehat{\pi}^{t2}_j(x;{\bf f^{t2}})\}_{j=1}^k$ and formulas \eqref{eq10}--\eqref{eq12} are valid. In these equalities, we replace $x$ with $\arccos x$, and then divide equalities \eqref{eq11} and \eqref{eq12} term by term by $\sqrt{1-x^2}$. Since the following identities hold on the segment $[-1,1]$ 
$$ 
f^{t2}_j(\arccos x)=\sqrt{1-x^2}f^{ch2}_j(x),\,\,\,j=1,\dots,k,
$$ 
then we obtain as a result:
$$
\widehat{Q}^t_m(\arccos x;{\bf f^{t2}})=\widehat{Q}_m(e^{i\arccos x};{\bf f^2}) \cdot \overline{\widehat{Q}_m(e^{i\arccos x};{\bf f^2})}=\sum_{l=0}^m{\widehat{q_l}T_l(x)},
$$
$$
\frac{1}{\sqrt{1-x^2}}\widehat{P}^t_j(\arccos x;{\bf f^{t2}})=\frac{1}{\sqrt{1-x^2}}\textrm{Im} \left \{ \widehat{P}_j(e^{i\arccos x};{\bf f^2}) \cdot \overline{\widehat{Q}_m(e^{i\arccos x};{\bf f^2})} \right \}=\sum_{l=0}^{n_j}{\widehat{p}^j_lU_l(x)},
$$
$$
f^{ch2}_j(x)-\frac{1}{\sqrt{1-x^2}}\widehat{\pi}^t_j(\arccos x;{\bf f^{t2}})=\sum_{l=n+m+1}^{\infty}{\widehat{d}^j_lU_l(x)}.
$$
This implies the existence of the Hermite--Chebyshev approximations $\{\widehat{\pi}^{ch2}_j(x;{\bf f^{ch2}})\}_{j=1}^k$ and the validity of equalities \eqref{eq17}, \eqref{eq18}. Theorem 4 is proved. Theorem 3 is proved similarly.

\vspace{0,2cm}
{\bf Remark 3.} For $k=1$, without describing the explicit form of the Pad\'e--Chebyshev approximations, theorem 3 is proved in \cite{Suetin1}. The assertions of theorem 4 are new and are of independent interest, including in the case when $k=1$. The fact is that nonlinear Pad\'e--Chebyshev approximations of the second kind, just like nonlinear Pad\'e--Chebyshev approximations of the first kind, are in demand in various applications \cite{Andr, Andr2, Kni, Ermoh}, as well as in scientific and technical calculations (for more details see \cite{Suetin1, Prox, Tee}).

We also note that for nonlinear Hermite--Chebyshev approximations of the second kind, the problem of finding existence conditions has not been researched to date. 

\renewcommand{\refname}{\center{\large\bf REFERENCES}}

\vspace{0,2cm}
\noindent {\it Aleksandr Pavlovich Starovoitov},
Dr. Phys.-Math. Sci., Prof.,
Francisk Scorina Gomel State University, Gomel, 246003 Belarus, \\
e-mail: svoitov@gsu.by \\

\noindent 
{\it Igor Viktorovich Kruglikov},
student,
Francisk Scorina Gomel State University, Gomel, 246003 Belarus, \\
e-mail: igor.v.kruglikov@gmail.com

\end{document}